\documentclass[reqno,11pt]{amsart}

\usepackage{amsthm, amsmath, amssymb, stmaryrd}
\usepackage[utf8]{inputenc}
\usepackage[T1]{fontenc}

\usepackage[hyphenbreaks]{breakurl}
\usepackage[hyphens]{url}
\usepackage{systeme}
\usepackage[shortlabels]{enumitem}
\usepackage[hidelinks]{hyperref}
\usepackage{microtype}

\usepackage{bm}
\usepackage[margin=1in]{geometry}

\usepackage[textsize=scriptsize,backgroundcolor=orange!5]{todonotes}

\usepackage[noabbrev,capitalize,sort,nameinlink]{cleveref}
\crefname{equation}{}{}
\numberwithin{equation}{section}

\usepackage{mathtools}
\newtheorem{theorem}{Theorem}[section]
\newtheorem{proposition}[theorem]{Proposition}
\newtheorem{lemma}[theorem]{Lemma}
\newtheorem{claim}[theorem]{Claim}
\newtheorem{corollary}[theorem]{Corollary}
\newtheorem{conjecture}[theorem]{Conjecture}

	\crefname{claim}{Claim}{Claims}

\theoremstyle{definition}
\newtheorem{definition}[theorem]{Definition}

\newtheorem{example}[theorem]{Example}

\theoremstyle{remark}

\newcommand{\abs}[1]{\left\lvert#1\right\rvert}

\DeclareMathOperator{\codim}{codim}

\DeclareMathOperator{\mult}{mult}   
\DeclareMathOperator{\vol}{vol}

\newcommand*{\eqdef}{\stackrel{\mbox{\normalfont\tiny def}}{=}}

\newcommand{\FF}{\mathbb{F}}
\newcommand{\RR}{\mathbb{R}}
\newcommand{\NN}{\mathbb{N}}
\newcommand{\ZZ}{\mathbb{Z}}

\newcommand{\cL}{\mathcal L}

\newcommand{\veps}{\varepsilon}
\newcommand{\Hasse}{\mathsf H} 

\newlength{\hght}

\makeatletter
\newcommand\thankssymb[1]{\textsuperscript{\@fnsymbol{#1}}}
\makeatother

\author[Ting-Wei Chao]{Ting-Wei Chao\thankssymb{1}}
\author[Hung-Hsun Hans Yu]{Hung-Hsun Hans Yu\thankssymb{2}}

\thanks{\thankssymb{1}Department of Mathematics, Massachusetts Institute of Technology, Cambridge, MA 02139\@. Email: {\tt twchao@mit.edu}}

\thanks{\thankssymb{2}Department of Mathematics, Princeton University, Princeton, NJ 08544\@. Email: {\tt hansonyu@princeton.edu}}

\title{Finite field Nikodym problem for spread line sets}
\begin{document}
\maketitle
\begin{abstract}
    A set of points $N\subseteq \FF_q^d$ is a Nikodym set if, for any $x\in \FF_q^d$, there is a line $\ell$ through $x$ such that $\ell\setminus\{x\}\subseteq N$. We conjecture that $\abs{N}=q^d-O_d(q^{d/(d-1)})$ and prove it under an extra algebraic assumption.
\end{abstract}
\section{Introduction}
Let $\FF_q$ be a finite field.
A \emph{Nikodym set} $N\subseteq \FF_q^d$ is a set such that for any $x\in \FF_q^d$, there is a line $\ell$ through $x$ so that $\ell\backslash\{x\}\subseteq N$.
A \emph{weak Nikodym set} $N\subseteq \FF_q^d$ is a set where we now only require for any $x\in \FF_q^d\backslash N$, there is a line $\ell$ throuogh $x$ with $\ell\backslash \{x\}\subseteq N$.
With those definitions, it is clear that a Nikodym set is always a weak Nikodym set.
Given a prime power $q$ and a positive integer $d$, the finite field Nikodym problem asks to determine the minimum possible size of a Nikodym set or a weak Nikodym set in $\FF_q^d$.

The finite field Nikodym problem is closely related to the \emph{finite field Kakeya problem}, asking to determine the smallest size of a \emph{Kakeya set} in $\FF_q^d$, namely a set containing a line in every direction.
The finite field Kakeya problem was first raised by Wolff \cite{Wol99} as a finite-field analog of the Kakeya set problem, a central problem in harmonic analysis that is still open in its full generality (see \cite{WZ25} nonetheless for a recent breakthrough by Wang and Zahl in three dimension).
In particular, Wolff conjectured that any Kakeya set in $\FF_q^d$ has size at least $c_dq^d$ for some constant $c_d>0$ only depending on $d$.
An elegant proof of this conjecture was given by Dvir \cite{Dvir09}, where he implicitly showed that union of dilates of a Kakeya set is a Nikodym set, and a Nikodym set has to be large (see \cite[Section 4.2.1]{Dvir10} where this is explicitly written out).
The proof uses the polynomial method, and it has inspired many later works in incidence geometry.
We refer the readers to \cite{Guth-book} and \cite{Sheffer-book} for more on the polynomial method.
We also refer the readers to \cite{HPVZ26} for other problems related to the finite field Nikodym problem.

Given the close connection between the finite field Nikodym problem and the finite field Kakeya problem, it is perhaps surprising how much more we know about the finite field Kakeya problem than the finite field Nikodym problem.
In fact, Bukh and the first author \cite{BC21} showed that Kakeya sets in $\FF_q^d$ have sizes at least $\left(\frac{1}{2^{d-1}}+o_{d;q\to\infty}(1)\right)q^d$, matching the construction by Dvir in \cite{SS08}; on the other hand, all known constructions of (weak) Nikodym sets have sizes $(1-o_{d;q\to\infty}(1))q^d$, and the best lower bound that holds for general $q,d$ is essentially the lower bound for the size of a Kakeya set in \cite{BC21} (see \cite{Tao25} for a deduction of this from \cite{BC21}).
It is still wide open whether the smallest (weak) Nikodym set has size $(1-o_{d;q\to\infty}(1))q^d$.
Currently, the strongest piece of evidence supporting that the upper bound is closer to the truth is a result by Guo, Kopparty and Sudan \cite{GKS13} that says the following: when restricted to finite fields with bounded characteristic, any Nikodym set in $\FF_q^d$ has size at least $q^d-q^{(1-\varepsilon)d}$ for some $\varepsilon>0$ depending on $d$ and the bound on the characteristic.

In this paper, we pose an even stronger conjecture and provide some initial evidence supporting our conjecture.

\begin{conjecture}
    For any $d\geq 3$, there exists a constant $C_d>0$ so that any weak Nikodym set in $\FF_q^d$ has size at least $q^d-C_dq^{d-\frac{1}{d}}$.
\end{conjecture}

We remark that if this is true, then the exponent $d-\frac{1}{d}$ would be tight by a very recent and independent work by Hunter, Pohoata, Verstraete and Zhang \cite[Theorem 1.8]{HPVZ26}.
We also remark that the analogous statement for $d=1$ is trivial, and the one for $d=2$ also holds (see \cite{CPZ25} and \cite{HPVZ26} for example).

To support our conjecture, we show that the conjecture is true if we assume that ``the lines associated with the weak Nikodym set'' are ``very generic''.
To be more specific, for any weak Nikodym set in $\FF^d_q$, we may associate with it a set of lines $\cL_N=\{\ell_x:x\in \FF^d_q\backslash N\}$ where $\ell_x\backslash \{x\}\subseteq N$ for any $x\in \FF^d_q\backslash N$.
Then we can prove our conjecture when the set of lines $\cL_N$ is ``algebraically spread''. As the definition is slightly technical, we will postpone its definition to \cref{sec:alg_spread}.

\begin{theorem}\label{thm:Nikodym}
    For any $d\in\NN$, there exists a constant $C_d>0$ so that the following holds.
    Let $N$ be a weak Nikodym set in $\FF^d_q$. 
    Let $\cL_N = \{\ell_x:x\in \FF^d_q\backslash N\}$ be a set of lines associated with the weak Nikodym set $N$.
    Let $\FF$ be any field extension of $\FF_q$ and suppose that the set of lines $\cL_N$ is algebraically spread in $\FF^d$. 
    Then $\abs{N}\geq q^d-C_dq^{d-\frac{1}{d}}$.
\end{theorem}
To give a concrete example of an algebraically spread set of lines, if $\abs{N} = q^d-a^{d-1}$ for some positive integer $a$, then the set of lines $\cL_N$ has size $a^{d-1}$.
Suppose that there is a hyperplane $H$ in $\FF^d$ and an invertible affine transformation $\phi:\FF^{d-1}\to H$ so that $H$ intersects lines in $\cL_N$ at $a^{d-1}$ distinct points, and those points are precisely $\phi(A_1\times \cdots \times A_{d-1})$ for some $A_1,\ldots, A_{d-1}\subseteq \FF$ each of size $a$.
Then we will see that $\cL_N$ is algebraically spread in $\FF^d$, and so $\abs{N}\geq q^d-C_dq^{d-\frac{1}{d}}$ in this particular case.
\\
\\
\noindent{\bf Structure of the paper.} We start with a brief review of Hasse derivatives and multiplicities in \cref{sec:prelim}. Then introduce the notion of algebraic spreadness in \cref{sec:alg_spread}. Finally, we prove our main theorem (\cref{thm:Nikodym}) in \cref{sec:main}.
\section{Preliminary}\label{sec:prelim}
We briefly recall the definitions and some properties of Hasse derivatives and multiplicities. For a more detailed introduction, see \cite{DKSS13}.
\begin{definition}[Hasse derivatives]
    For any $f\in\FF[x_1,\dots,x_d]$, its Hasse derivatives $\Hasse^\alpha f(x),\,\alpha\in\ZZ_{\geq 0}^d$ are given by
    \[f(x+y)=\sum_{\alpha\in\ZZ_{\geq 0}^d}\Hasse^\alpha f(x)y^\alpha.\]
\end{definition}
\begin{definition}[Multiplicity]
    Let $f\in\FF[x_1,\dots,x_d]$ and $p\in \FF^d$. The multiplicity $\mult(f,p)$ of $f$ at $p$ is the largest $n$ such that $\Hasse^\alpha f(p)=0$ holds for all $\abs{\alpha}<n$. If this holds for every $n\in \NN$, then we define $\mult(f,p)=\infty$.
\end{definition}
\begin{definition}[Restriction of a polynomial]
     Let $f\in\FF[x_1,\dots,x_d]$. Suppose that $F$ is a $r$-dimensional affine subspace in $\FF^d$. We denote by $f|_F$ the polynomial $f(a+b_1t_1+\dots+b_rt_r)$, where $\{a+b_1t_1+\dots+b_rt_r\mid t_1,\dots,t_r\in\FF\}$ is a fixed parameterization of $F$. 
\end{definition}
\begin{definition}[Multiplicity on a line]
    Let $f\in\FF[x_1,\dots,x_d]$. Suppose that $\ell$ is a line in $\FF^d$ and $p\in\ell$.
    The multiplicity $\mult(f|_{\ell},p)$ is the largest $n$ such that $\Hasse^k (f|_{\ell})=0$ at $p$ holds for all $k<n$. If this holds for every $n\in \NN$, then we define $\mult(f|_{\ell},p)=\infty$.
\end{definition}
\begin{proposition}
    Multiplicities are invariant under affine transformation.
\end{proposition}
\begin{proposition}
    Let $f\in\FF[x_1,\dots,x_d]$. Suppose that $\ell$ is a line in $\FF^d$ and $p_1,\dots,p_r\in\ell$. If 
    \[\sum_{i=1}^r \mult(f|_\ell,p_i)>\deg f,\]
    then $f|_\ell$ must be the zero polynomial.
\end{proposition}

\section{Algebraic spreadness}\label{sec:alg_spread}
Below, we let $\FF$ be any field extension of $\FF_q$.
We will eventually take $\FF$ to be any uncountable field extension of $\FF_q$.
In this section, we will define what it means for a set of lines to be algebraically spread.
We begin by defining what it means for a set of points to be algebraically spread.

\begin{definition}[Algebraic spreadness for points]\label{def:alg_spread_point}
    We say that a set of points $\{p_1,\dots,p_k\}$ on an $r$-dimensional flat $F$ in $\FF^d$ is \emph{algebraically spread} if it satisfy the following property. For any $\veps>0$, there exists an integer $n_0$ such that the following holds for any positive integer $n>n_0$. For any polynomial $f$ defined on $F$ with degree less than $(1-\veps)k^{1/r}n$, if $\mult(f,p_i)\geq n$ holds for all $i\in [k]$, then $f$ must be the zero polynomial.
\end{definition}

\begin{example}
    In $F\cong \FF^r$, consider any $A_1\times \cdots\times A_r\subseteq \FF^r$ where $A_1,\ldots, A_r$ each has size $a$.
    Then $k=a^r$.
    Moreover, by the generalized Schwartz--Zippel lemma (see \cite[Lemma 2.7]{DKSS13}), we know that any polynomial $f$ with degree less than $an$ and $\mult(f,p)\geq n$ for all $p\in A_1\times \cdots\times A_r$ must be the zero polynomial. Therefore, the point set $A_1\times \cdots\times A_r$ is algebraically spread.
\end{example}

We briefly remark that we shall expect generic points to satisfy algebraic spreadness for the following reason: locally at each $p_i$, the assumption that $\mult(f,p_i)\geq n$ can be viewed as $\binom{n+r-1}{r}$ linear conditions, and there are $\binom{(1-\varepsilon)k^{1/r}n+r-1}{r}$ degrees of freedom for $f$.
Therefore if the linear conditions are almost linearly independent, then since the total number of linear conditions $k\binom{n+r-1}{r}$ is more than the dimension $\binom{(1-\varepsilon)k^{1/r}n+r-1}{r}$, we may conclude that $f$ has to be zero.
However, whether this intuition holds true for sufficiently large $k$ is widely open: in fact even when $r=2$ and $k>9$ is not a square, this is closely related to a well-known conjecture by Nagata \cite{Nagata59}. 

Before we move on to lines, let us prove some basic properties of algebraic spreadness for points.
Below, a \emph{very generic} object (for example set of points or set of lines) means one that is chosen from a countable union of proper Zariski-closed set in the space of object.

\begin{proposition}\label{prop:very-generic-point}
Let $\{p_1,\ldots, p_k\}\subseteq \FF^r$ be an algebraically spread set of points.
Then it remains algebraically spread in $(\FF')^r$ where $\FF'$ is any field extension of $\FF$.
Moreover, a very generic set of $k$ points in $\FF^r$ is algebraically spread.
    
\end{proposition}
\begin{proof}
    By definition, we may find an infinite sequence $D_1,D_2,D_3,\ldots$ of positive integers so that $\lim_{n\to\infty}D_n/n = k^{1/r}$, and that for any $n\in\NN$ if a polynomial $f$ has degree less than $D_n$ and multiplicity at least $n$ at $p_1,\ldots, p_k$ then $f$ must be zero.

    Let us first show that this does not depend on the base field.
    To see this, simply note that the above condition is equivalent to that the matrix
    \[M_n\eqdef \left((\Hasse^\alpha x^\beta)\mid_{x=p_i}\right)\]
    is full rank where columns are indexed by $\alpha\in \ZZ^r_{\geq 0}$ with $\abs{\alpha}<D_n$ and rows are indexed by $(\beta,i)\in \ZZ^r_{\geq 0}\times [k]$ with $\abs{\beta}<n$.
    Therefore this is field-independent.

    Moreover, we see that when viewing $(p_1,\ldots, p_k)$ as a variable, the set $S_n$ of $(p_1,\ldots, p_k)\in (\FF^r)^k$ where $M_n$ is of full-rank is non-empty and Zariski-open.
    Therefore, any set of points in $\bigcap_{n=1}^{\infty}S_n$ is algebraically spread.
    As $\bigcap_{n=1}^{\infty}S_n$ is very generic, the statement follows.
\end{proof}

We can now lift the definition of algebraic spreadness to lines.
\begin{definition}[Algebraic spreadness for lines]\label{def:alg_spread_line}
    We say that a set of lines $\cL=\{\ell_1,\dots,\ell_k\}$ in $\FF^d$ is \emph{algebraically spread} if there exists a hyperplane $H$ such that $\cL\cap H\eqdef\{\ell_i\cap H\mid i\in [k]\}$ is an algebraically spread point set on $H$ where $\ell_i\cap H\neq \ell_j\cap H$ for any $i\neq j\in [k]$.
\end{definition}

As a generic hyperplane $H$ intersects the lines in $\cL$ at distinct points, and those points depend polynomially on $H$, we immediately have the following corollary from \cref{prop:very-generic-point}.

\begin{corollary}\label{cor:very-generic-hyperplane}
    Let $\cL$ be a set of lines in $\FF^d$ that is algebraically spread, and let $\FF'$ be any field extension of $\FF$.
    Then $\cL$ is algebraically generic in $(\FF')^d$ as well.
    Moreover, for a very generic hyperplane $H$ in $\FF^d$, we have that $\cL\cap H$ is a set of $\abs{\cL}$ disjoint points that is algebraically spread.
\end{corollary}

\section{Main proof}\label{sec:main}
In this section, we give a proof of \cref{thm:Nikodym} by proving a slightly stronger statement as follows.
\begin{theorem}\label{thm:NL}
    Let $\cL$ be a set of lines in $\FF^d$ that is algebraically spread. Suppose $H$ is a hyperplane such that $\cL\cap H$ is an algebraically spread point set on $H$ of size $\abs{\cL}$. Let $N_\ell$ be a set of $q-1$ points on $\ell$ for each $\ell\in\cL$ such that $N_{\ell}\cap H=\emptyset$, and let $N=\cup_{\ell\in\cL}N_\ell$. Assume $\abs{N}+\abs{\cL}\leq q^d$. Then $\abs{\cL}=O_d(q^{d-\frac{1}{d}})$.
\end{theorem}
We first show that \cref{thm:NL} implies \cref{thm:Nikodym}.
\begin{proof}[Proof of \cref{thm:Nikodym} assuming \cref{thm:NL}]
    Let $N$ be a weak finite field Nikodym set in $\FF_q^d$ such that $\cL_N=\{\ell_x\mid x\in \FF_q^d\setminus N\}$ is algebraically spread in $\FF$.
    By \cref{cor:very-generic-hyperplane}, we may extend $\FF$ if necessarily to assume that $\FF$ is uncountable.
    As $\FF$ is uncountable, we may also take $H$ to avoid $\FF_q^d$ while maintaining that $\cL_N\cap H\eqdef\{\ell\cap H\mid \ell\in\cL_N\}$ is an algebraically spread point set on $H$ of size $\abs{\cL_N}$.
    By the definition of weak finite field Nikodym sets, we know that $N\cap \ell$ is a set of $q-1$ points for all $\ell\in\cL_N$. Thus, we may apply \cref{thm:NL} with  the point set $N'=\cup_{\ell\in\cL_N}N\cap\ell$ and the line set $\cL_N$. Since $\abs{N'}+\abs{\cL_N}\leq \abs{N}+\abs{\FF_q^d\setminus N}=q^d$, we have $\abs{\FF_q^d\setminus N}=\abs{\cL}=O_d(q^{d-\frac{1}{d}})$, and the theorem follows.
\end{proof}
Now it remains to prove \cref{thm:NL}.
\begin{proof}[Proof of \cref{thm:NL}]
    Note that the assumptions in the theorem are invariant under affine transformations. Thus, we may apply an affine transformation and assume without loss of generality that $H=\{x_d=0\}$.
    Since $\cL\cap H$ is a set of points, we know that all lines in $\cL$ are not parallel to $H$. We call such lines \emph{non-horizontal}.

    We begin with a vanishing lemma which is analogous to the one in \cite{BC21}. To do so, we need the following definition of vanishing introduced in \cite{BC21}.
    \begin{definition}
        Let $\ell$ be a non-horizontal line and $p\in \ell$. We say that a polynomial $f\in\FF[x_1,\dots,x_d]$ vanishes to order $(u,v;\gamma)$ on $p$ along $\ell$ with if $\mult((\Hasse^{(\alpha,0)}f)|_\ell,p)\geq v-\frac{\abs{\alpha}}{\gamma}$ holds for all $\alpha\in\ZZ_{\geq 0}^{d-1}$ with $\abs{\alpha}<u$.
    \end{definition}
    Let $V(r,s)$ be the space of polynomial spanned by monomials $x^{\alpha}$, where $\alpha=(\alpha_1,\dots,\alpha_d)\in \ZZ_{\geq 0}^d$ with $\alpha_1+\dots+\alpha_{d-1}<r$ and $\abs{\alpha}<s$.
    \begin{lemma}\label{lemma:vanish}
        Let $c>0$, $\veps>0$, and $n\in\NN$ sufficiently large.
        If $f\in V\left((1-\veps)\abs{\cL}^{1/(d-1)}n,(q-1)cn\right)$ vanishes to order $(n,cn;q-1)$ on $p$ along $\ell$ for every $p\in N_\ell$ and $\ell\in\cL$, then $f$ must be the zero polynomial.
    \end{lemma}
    \begin{proof}
        We will prove by inducting on $k$ that $f$ is divisible by $x_d^k$. It is true for $k=0$. Assume the statement is true for $k$, now we show it for $k+1$. Set $g=f/x_d^k$ and note that we still have $g\in V\left((1-\veps)\abs{\cL}^{1/(d-1)}n,(q-1)cn\right)$. Moreover, we know that $g$ still vanishes to order $(n,cn;q-1)$ on $p$ along $\ell$ with for every $p\in N_\ell$ and $\ell\in\cL$ since $x_d^k$ is not zero on any point in $N_{\ell}$.

        Now, it is enough to show that $g$ is divisible by $x_d$. For any $\alpha\in\ZZ_{\geq 0}^{d-1}$ with $\abs{\alpha}<n$, we know that $(\Hasse^{(\alpha,0)}g)|_\ell$ has multiplicity at least $cn-\frac{\abs{\alpha}}{q-1}$ on each point in $N_{\ell}$. Since there are $q-1$ points in $N_{\ell}$, the total multiplicities is at least $(q-1)cn-\abs{\alpha}$. On the other hand, we have 
        \[\deg(\Hasse^{(\alpha,0)}g)|_\ell\leq \deg\Hasse^{(\alpha,0)}g\leq \deg g-\abs{\alpha}<(q-1)cn-\abs{\alpha},\]
        and hence $(\Hasse^{(\alpha,0)}g)|_\ell$ must be the zero polynomial. In particular, we know that $\Hasse^{(\alpha,0)}g(p_{\ell})=0$ for any $\ell\in\cL$, where $p_{\ell}=\ell\cap \{x_d=0\}$. That is, $\mult(g|_H,p_\ell)\geq n$ for all $\ell\in\cL$. From the definition of algebraic spreadness of $\cL$, we can conclude that $g|_H$ must be the zero polynomial. Thus, $g$ is divisible by $x_d$.
    \end{proof}

    Next, we proceed with a dimension counting analogous to the one in \cite{BC21}. For each $p\in N$, let $m_p$ be the number $\ell\in\cL$ such that $p\in N_\ell$. Then we have the following inequality.
    \begin{lemma}\label{lemma:dim_counting}
    We have
    \[\sum_{p\in N}m_p\left(m_p^{1/(d-1)}-1\right)\leq \abs{\cL}\left(\abs{\cL}^{1/(d-1)}-1\right).\]
    \end{lemma}
    \begin{proof}
    We set $c$ to be the maximum among $\abs{\cL}^{1/(d-1)}/(q-1)$ and $(m_p^{1/(d-1)})_{p\in N}$, and define the vector space
    \[C_p=\{f\in\FF[x_1,\dots,x_d]\mid f\text{ vanishes to order $(n,cn;q-1)$ on $p$ along $\ell$ for every $\ell$ with $p\in N_{\ell}$}\}.\]
    It follows from \cref{lemma:vanish} that
    \begin{align}\label{eq:codim}
        \dim V\leq \sum_{p\in N}\codim_{V} \left(C_p\cap V\right)\leq \sum_{p\in N}\codim_{\FF[x_1,\dots,x_d]}C_p,
    \end{align}
    where $V=V\left((1-\veps)\abs{\cL}^{1/(d-1)}n,(q-1)cn\right)$. 
    
    In the remaining part of the proof of this lemma, the $o(1)$ factor may depend on $q,d,\cL,N$ but not on $n$, and goes to $0$ as $n$ goes to infinity.  
    
    We first compute the left hand side of \cref{eq:codim}.
    \begin{claim}\label{claim:dimV}
         We have
    \begin{align*}
        \dim V=\left(1-o(1)\right)\left(-\frac{d-1}{d!}\abs{\cL}^{d/(d-1)}+\frac{1}{(d-1)!}(q-1)c\abs{\cL}\right)n^d.
    \end{align*}
    \end{claim}
    \begin{proof}
    Note that as $n$ goes to infinity, the parameter $\veps>0$ can be chosen arbitrarily small, so we have
    \[\dim V=\left(1-o(1)\right)\dim V\left(\abs{\cL}^{1/(d-1)}n,(q-1)cn\right).\] 
    If we set 
    \[S(u,v)=\{a=(a_1,\dots,a_d)\in \RR_{\geq 0}^d\mid a_1+\dots+a_{d-1}\leq u,\abs{a}\leq v\},\]
    then we have
    \[\dim V(un,vn)=(1-o(1))\vol(S(u,v))n^d.\]
    If $v\geq u$, we can compute the volume of $S(u,v)$ directly and get
    \[\vol(S(u,v))=\frac{u^d}{d!}+(v-u)\frac{u^{d-1}}{(d-1)!}=-\frac{d-1}{d!}u^d+\frac{1}{(d-1)!}vu^{d-1}.\]
    Therefore, by taking $u=\abs{\cL}^{1/(d-1)}$ and $v=(q-1)c$, and note that $v\geq u$ since $c\geq \abs{\cL}^{1/(d-1)}/(q-1)$, the claim follows.
    \end{proof}

    Next, we compute the right hand side of \cref{eq:codim}.
    \begin{claim}\label{claim:Cp}
        We have
    \begin{align*}
        \codim_{\FF[x_1,\dots,x_d]}C_p\leq \left(1+o(1)\right)\left(-\frac{d-1}{d!}m_p\left(m_p^{1/(d-1)}-1\right)+m_p\left(c\frac{1}{(d-1)!}-\frac{1}{q-1}\frac{d-1}{d!}\right)\right)n^d.
    \end{align*}
    \end{claim}
    \begin{proof}
        For each $p\in N$ and each $\ell\in \cL$ such that $p\in N_\ell$, we define
    \[C^-_p=\{f\in \FF[x_1,\dots,x_d]\mid\mult(f,p)\geq m_p^{1/(d-1)}n\}\]
    and
    \[C_{p,\ell}=\{f\in C^-_p\mid f\text{ vanishes to order $(n,cn;q-1)$ on $p$ along $\ell$}\}.\]

    Note that 
    \begin{align}\label{eq:Cp_leq_Cp-+Cpl}
        \codim_{\FF[x_1,\dots,x_d]}C_p\leq\codim_{\FF[x_1,\dots,x_d]}C^-_p+\sum_{\ell:p\in N_{\ell}}\codim_{C^-_p}C_{p,\ell}.
    \end{align}
    From a direct calculation, we have
    \begin{align}\label{eq:Cp-}
        \codim_{\FF[x_1,\dots,x_d]}C^-_p=\left(1+o(1)\right)\frac{m_p^{d/(d-1)}n^d}{d!}.
    \end{align}
    To compute $\codim_{C^-_p}C_{p,\ell}$, we may assume without loss of generality that $\ell$ is parallel to the $x_d$-axis since the codimension is invariant under affine transformations. In this case, we know that 
    \[C_{p,\ell}\subseteq\Bigl\{f\in C^-_p \,\Bigm|\,\Hasse^{(\alpha,\beta)}f(p)=0
\begin{aligned}[t]
&\text{ for all } \alpha\in\ZZ_{\geq 0}^{d-1},\beta\in\ZZ_{\geq 0}\\
&\text{ such that }\abs{\alpha}+\beta\geq m_p^{1/(d-1)}n,\abs{\alpha}<n,\beta< cn-\frac{\abs{\alpha}}{q-1}
\Bigr\}.
\end{aligned}
\]
    Note the each $\Hasse^{(\alpha,\beta)}f(p)=0$ is a linear equation to coefficients of $f$. Thus, the codimension
    \begin{align*}
        \codim_{C^-_p}C_{p,\ell}\leq& \#\left\{(\alpha,\beta)\in\ZZ_{\geq 0}^{d-1}\times  \ZZ_{\geq 0}\,\middle|\, \abs{\alpha}+\beta\geq m_p^{1/(d-1)}n,\abs{\alpha}<n,\beta< cn-\frac{\abs{\alpha}}{q-1}\right\}\\
        =&(1+o(1))n^d\vol(T),
    \end{align*}
    where
    \[T=\left\{(a,b)\in\RR_{\geq 0}^{d-1}\times\RR_{\geq 0}\,\middle|\,\abs{a}+b\geq m_p^{1/(d-1)},\abs{a}<1,b< c-\frac{\abs{a}}{q-1}\right\}.\]
    We write
    \[T_1 = \left\{(a,b)\in\RR^{d-1}_{\geq 0}\times \RR_{\geq 0}\,\middle|\,\abs{a}+b\in \left[m_p^{1/(d-1)},c\right], \abs{a}<1\right\}\]
    and
    \[T_2 = \left\{(a,b)\in\RR^{d-1}_{\geq 0}\times \RR_{\geq 0}\,\middle|\,\abs{a}\in \left[c-b,(q-1)(c-b)\right), \abs{a}<1 \right\}.\]
    It is clear that $\vol(T) = \vol(T_1)+\vol(T_2)$.
    By a direct computation, we have
    \[\vol(T_1) = \left(c-m_p^{1/(d-1)}\right)\frac{1}{(d-1)!}\]
    (where we use $c\geq m_p^{1/(d-1)}\geq 1$) and
    \[\vol(T_2) = (d-1)\left(1-\frac{1}{q-1}\right)\frac{1}{d!}, \]
    giving
    \[\vol(T)=\left(c-m_p^{1/(d-1)}\right)\frac{1}{(d-1)!}+(d-1)\left(1-\frac{1}{q-1}\right)\frac{1}{d!}.\]
    Together with \cref{eq:Cp_leq_Cp-+Cpl,eq:Cp-}, we can conclude that
    \begin{align*}
        &\codim_{\FF[x_1,\dots,x_d]}C_p\\
        \leq&\left(1+o(1)\right)\left(\frac{m_p^{d/(d-1)}}{d!}+m_p\left((c-m_p^{1/(d-1)})\frac{1}{(d-1)!}+(d-1)\left(1-\frac{1}{q-1}\right)\frac{1}{d!}\right)\right)n^d,
    \end{align*}
    and the claim follows from rearranging the right hand side.
    \end{proof}
    By combining \cref{eq:codim,claim:dimV,claim:Cp}, using the fact that $\sum_{p\in N}m_p=(q-1)\abs{\cL}$, and taking $n$ to infinity, and we get
    \begin{align*}
        &-\frac{d-1}{d!}\abs{\cL}^{d/(d-1)}+\frac{1}{(d-1)!}(q-1)c\abs{\cL}\\
        \leq &-\sum_{p\in N}\frac{d-1}{d!}m_p\left(m_p^{1/(d-1)}-1\right)+(q-1)\abs{\cL}\left(c\frac{1}{(d-1)!}-\frac{1}{q-1}\frac{d-1}{d!}\right).
    \end{align*}
    After rearranging, we have
    \[\sum_{p\in N}m_p(m_p^{1/(d-1)}-1)\leq \abs{\cL}(\abs{\cL}^{1/(d-1)}-1).\qedhere\]
    \end{proof}

    We can rewrite \cref{lemma:dim_counting} using $\sum_{p\in N}m_p=(q-1)\abs{\cL}$ and get
    \[\sum_{p\in N}m_p^{d/(d-1)}\leq \abs{\cL}^{d/(d-1)}+(q-2)\abs{\cL}.\]
    Let $x=\abs{\cL}$. From the assumption, we know that $\abs{N}\leq q^d-x$. We may view $\sum_{p\in N}m_p^{d/(d-1)}$ as a summation of $q^d-x$ terms by adding $0^{d/(d-1)}$. By the concavity of $\star^{d/(d-1)}$, we have 
    \begin{align*}
        x^{d/(d-1)}+(q-2)x\geq& \sum_{p\in N}m_p^{d/(d-1)}\\
        \geq& \left(q^d-x\right)\left(\frac{(q-1)x}{q^d-x}\right)^{d/(d-1)}=\left(q^d-x\right)^{-1/(d-1)}(q-1)^{d/(d-1)}x^{d/(d-1)}.
    \end{align*}
    Dividing both sides by $x^{d/(d-1)}$ and rewriting the right hand side, we get
    \[1+(q-2)x^{-1/(d-1)}\geq\left(1-\frac{x}{q^d}\right)^{-1/(d-1)} \left(1-\frac{1}{q}\right)^{d/(d-1)}.\]
    Using Bernoulli's inequality, we have
    \[\left(1-\frac{x}{q^d}\right)^{1/(d-1)}\leq 1-\frac{1}{d-1}\cdot\frac{x}{q^d}\]
    since $-\frac{x}{q^d}\geq -1$ and $0\leq \frac{1}{d-1}\leq 1$, and also
    \[\left(1-\frac{1}{q}\right)^{d/(d-1)}\geq 1-\frac{d}{d-1}\cdot\frac{1}{q}\]
    since $-\frac{1}{q}\geq -1$ and $\frac{d}{d-1}\geq 1$.
    We may also assume that $x=\omega_{d;q\rightarrow\infty}(q^{d-1})$, otherwise the theorem holds vacuously.
    Therefore, we have
    \[1+(q-2)x^{-1/(d-1)}\geq \left(1-\frac{1}{d-1}\cdot\frac{x}{q^d}\right)^{-1}\left(1-\frac{d}{d-1}\cdot\frac{1}{q}\right)\geq 1+\frac{1}{d-1}\cdot\frac{x}{q^d}-\frac{d}{d-1}\cdot\frac{1}{q},\]
    where the second inequality follows from the fact that $(1-A)^{-1}(1-B)\geq 1+A-B$ holds for all $A\geq B\geq 0$.
     As $x = \omega_{d;q\to\infty}(q^{d-1})$, we can simplify the inequality above and get
    \[(q-2)x^{-1/(d-1)}\geq (1-o_{d;q\rightarrow\infty}(1))\frac{1}{d-1}\cdot\frac{x}{q^d}.\]
    After rearranging, we can conclude that
    \[x\leq (1+o_{d;q\rightarrow\infty}(1))(d-1)^{(d-1)/d}q^{d-1}(q-2)^{(d-1)/d}=O(q^{d-1/d}).\qedhere\]
    
\end{proof}

\bibliographystyle{plain}
\bibliography{bib.bib}
\end{document}